\documentclass[10pt]{article}

\usepackage{amsthm,amsmath,amssymb,amsbsy,amsfonts,graphics,epsfig}
\newcounter{geqncount} %
{\setcounter{equation}{\value{geqncount}}}

\newtheorem{prop}{Proposition}[section]
\newtheorem{theo}{Theorem}[section]

\newtheorem{lemma}{Lemma}[section]

\newcommand{\beq}{\begin{equation}} \newcommand{\eeq}{\end{equation}}
\newcommand{\beqr}{\begin{eqnarray}}
\newcommand{\eeqr}{\end{eqnarray}}
\newcommand{\beqrn}{\begin{eqnarray*}}
\newcommand{\eeqrn}{\end{eqnarray*}} \newcommand{\brr}{\begin{array}}
\newcommand{\err}{\end{array}}
\newcommand{\bseq}{\begin{subequations}}
\newcommand{\eseq}{\end{subequations}}
\newcommand{\bef}{\begin{figure}} \newcommand{\eef}{\end{figure}}
\newcommand{\bec}{\begin{center}} \newcommand{\eec}{\end{center}}

\newcommand{\R}{{\mathbb R}}

\title{Interlacing and non-orthogonality of spectral polynomials for the Lam\'e operator}

\author{A. Bourget, T. McMillen and A. Vargas}

\begin{document}

\maketitle

\bibliographystyle{plain}

\hrule

\begin{abstract}

\noindent Polynomial solutions to the Heine-Stieltjes equation, \textit{the Stieltjes polynomials}, and the
associated \textit{Van Vleck polynomials} have been studied since the 1830's in various contexts including the
solution of the Laplace equation on an ellipsoid.  Recently there has been renewed interest in the distribution
of the zeros of Van Vleck polynomials as the degree of the corresponding Stieltjes polynomials increases.  In
this paper we show that the zeros of Van Vleck polynomials corresponding to Stieltjes polynomials of successive
degrees interlace.   We also show that the spectral polynomials formed from the Van Vleck zeros
 are not orthogonal with respect to any measure.  This furnishes a counterexample, coming from a
 second order differential equation, to the well known theorem that the zeros of orthogonal polynomials interlace.

\vskip 10pt \hrule \vskip 10pt

\end{abstract}

\section{Introduction}
\label{s1}

Let $\alpha_1,\ldots,\alpha_n$ be any $n$ distinct complex numbers, and let $\rho_1,\ldots,\rho_n$ be positive
numbers. The \textit{generalized Lam\'e equation} is the second order ODE given by
\begin{equation} \label{Lame}
  \prod_{j=1}^n (z-\alpha_j) \, \phi''(z) + 2\sum_{j=1}^n  \rho_j \prod_{i \neq j}
  (z-\alpha_i) \, \phi'(z) = V(z) \, \phi(z).
\end{equation}
According to a result of Heine \cite{MR0204726}, there exist at most $\sigma(n,k)=\frac{(n+k-2)!}{(n-2)! \, k!}$
polynomials $V$ of degree $n-2$ for which \eqref{Lame} has a polynomial solution $\phi$ of degree $k$. These
polynomial solutions are often called \textit{Stieltjes polynomials}, and the corresponding polynomials $V$ are
known as \textit{Van Vleck polynomials}.

The equation \eqref{Lame} was studied by Lam\'e in the 1830's in the special case $n=3$, $\rho_i=1/2$,
$\alpha_1+\alpha_2+\alpha_3 = 0$ in connection with the separation of variables in the Laplace equation using
elliptical coordinates \cite[Ch. 23]{whwa96}. The equation has since found other applications in studies as
diverse as electrostatics and the quantum asymmetric top.

For $\alpha_1,...,\alpha_n$ real, Stieltjes \cite{MR1554669} showed that the location of the zeros of the
Stieltjes polynomials are completely characterized by their distribution in the subintervals
$(\alpha_1,\alpha_2),\ldots,(\alpha_{n-1},\alpha_n)$. Similar results for the zeros of the Van Vleck polynomials
were also obtained by Shah \cite{MR0230954}. Much is known about the properties of Van Vleck polynomials for a
fixed degree of the corresponding Stieltjes polynomial (see, e.g. \cite{vo04} for recent results), but there are
few results relating the Van Vleck zeros that correspond to Stieltjes polynomials of different degrees.
Recently there has been interest in the distribution of the zeros of Van Vleck polynomials as the degree of the
corresponding Stieltjes polynomials tends toward infinity \cite{bosh08}.

In this paper we consider the case of three $\alpha_i$'s on the real line and first degree Van Vleck polynomials.  In this case, let $\alpha_1<\alpha_2<\alpha_3$ and $\rho_1, \rho_2, \rho_3 > 0$   be real numbers and define
\[A(x) = \prod_{j=1}^3\left(x-\alpha_j\right), \;\;\;
B(x) =  \sum_{j=1}^3 2\rho_j \prod_{i\neq j}  (x-\alpha_i).
\]
Then the Lam\'e equation is
\begin{equation}
A(x)\phi'' + B(x)\phi' = \mu (x-\nu) \phi.
\label{lame1}
\end{equation}
By the Van Vleck zeros of order $k$ we mean the set of all $\nu$'s such that \eqref{lame1} has a polynomial solution of degree $k$.  In this case where the $\alpha_i$'s are real, Heine's result is exact, and Van Vleck showed that the $k+1$ Van Vleck zeros of order $k$ are distinct and lie in the interval $(\alpha_1,\alpha_3)$ \cite{vv98}.

This paper has two main results.  In \S2
we show that the Van Vleck zeros of successive orders interlace.  That is, if the Van Vleck zeros of order $k$ are written in increasing order as
 $\nu_1^{(k)} < \nu_2^{(k)} < \cdots < \nu_{k+1}^{(k)}$, then
\begin{equation}
\alpha_1 < \nu_{1}^{(k+1)} < \nu_1^{(k)} < \nu_{2}^{(k+1)} < \nu_2^{(k)} < \cdots < \nu_{k+1}^{(k)} < \nu_{k+2}^{(k+1)} < \alpha_3.
\label{inequality}
\end{equation}
The proof of this result will be carried out in two steps.  First we will show that the Van Vleck zeros of order $k$ and $k+1$ are distinct.  Then we will show that the interlacing property \eqref{inequality} holds for a special set of $\alpha_i, \rho_i$.  Since the Van Vleck zeros are continuous functions of these parameters, the interlacing property must hold in general.

Given the interlacing property and the well-known properties of orthogonal polynomials, it is natural to ask whether the polynomials formed from the Van Vleck zeros,
\[ \prod_{i=1}^{k+1} \left( x - \nu_i^{(k)}\right),
\]
are orthogonal with respect to some measure.  These polynomials are often referred to as \textit{spectral
polynomials}, as the Van Vleck zeros are often interpreted as an energy.  Their study goes back as far as
Hermite, and recent results are found in \cite{ma08} and the references therein.  In \cite{grve08, agbo07} the
spectral polynomials are studied in conjunction with the quantum Euler top.  The authors of \cite{grve08}, in
particular, relate the coefficients of the spectral polynomials to Bernoulli polynomials. In \S3 we show that
the spectral polynomials are not orthogonal with respect to any measure.  We conclude with a few remarks and a
conjecture.

\section{Van Vleck zeros of successive orders interlace}
\label{s2}

Our first step is to prove a lemma showing that the Van Vleck zeros of order $k$ and $k+1$ are distinct.  We follow a Sturm comparison type argument to argue that if there were a Van Vleck zero of order $k$ and $k+1$ in common, the Stieltjes polynomials corresponding to this common zero would have interlacing zeros, leading to a contradiction.  Before we proceed to the lemma we collect some facts about Stieltjes polynomials that will be used throughout.

The first is a result due to Stieltjes \cite{MR1554669} (see also \cite{po12}), that the zeros of every polynomial solution $S(z)$ of \eqref{Lame} lie in the smallest convex polygon containing $\alpha_1,\dots,\alpha_n$.  This follows from the fact that if $z_1,\dots, z_l$ are the zeros of a polynomial $S$ satisfying \eqref{Lame}, and $z_r$ is not an $\alpha_i$, then
\[ \sum_{j\neq r} \frac{1}{z_r - z_j} + \sum_{j} \frac{\rho_j}{z_r - \alpha_j} = 0.
\]
Appealing to the Gauss-Lucas Theorem, $z_r$ must lie in the smallest convex polygon containing $z_1,\dots, z_{r-1}, z_{r+1}, \dots, z_l, \alpha_1, \dots, \alpha_n$.  Since this is true for each $r$, the result follows.
Additionally, every zero $z_i$ of $S$ is simple unless $z_i$ corresponds to an $\alpha_i$, otherwise, if $S'(z_r)=0$ then repeated differentiation of \eqref{Lame} would show that $S$ is identically zero.
Szeg{\H{o}} \cite{sz75} showed that when the $\alpha_i$'s are real,  as for the equation \eqref{lame1} we are considering, Stieltjes polynomials cannot have zeros at any of the $\alpha_i$'s.  Thus, when the $\alpha_1<\alpha_2<\alpha_3$ are real, the zeros of Stieltjes polynomials are simple and real, and lie in  $(\alpha_1, \alpha_3)\setminus\{\alpha_2\}$.

Secondly, we will use two results due to Shah \cite{MR0230954}.  For one,
 no Van Vleck zero is a zero of the corresponding Stieltjes polynomial.   And, between any zero of a Stieltjes polynomial and a zero of a corresponding Van Vleck polynomial there is either a zero of the derivative of the Stieltjes polynomial or a singular point $\alpha_i$.

Finally, we note a result originally proved by Van Vleck \cite{vv98} in the real case, and extended to the complex case by Marden \cite{MR1501624} that every zero of the Van Vleck polynomials $V(z)$ in \eqref{Lame} lies in the smallest convex polygon containing $\alpha_1,\dots, \alpha_n$.  For the equation under consideration here \eqref{lame1}, every $\nu$ for which \eqref{lame1} admits a polynomials solution lies in $(\alpha_1, \alpha_3)$.

We now proceed to the lemma, which is completely general for all $\rho_i > 0$ and $\alpha_1<\alpha_2<\alpha_3$.  The equation \eqref{lame1} is invariant under affine transformations $x\mapsto ax+b$, so for convenience we assume that  $\alpha_1< 0 = \alpha_2 < \alpha_3$.
\begin{lemma}
Let $k$ be a positive integer.  Then no Van Vleck zero of order $k$ is a Van Vleck zero of order $k+1$.
\label{unequal}
\end{lemma}
\begin{proof}
Suppose there is a $\nu$ such that $S_k$ and $S_{k+1}$ are Stieltjes polynomials of degree $k$ and $k+1$, respectively, corresponding to Van Vleck polynomials with a zero at $\nu$.  Then $\nu\in(\alpha_1, \alpha_3)$, and since \eqref{lame1} is invariant under $x\mapsto -x$, we may assume that $\nu\in [0, \alpha_3)$.  Employing the results noted above, the zeros of $S_k$ and $S_{k+1}$ are all simple and in $(\alpha_1, \alpha_3)\setminus\{0\}$.  Since between any zero of $S_i$  $(i\in\{k,k+1\})$ and $\nu$ there is either a zero of $S_i'$ or 0, all of the zeros of $S_k$ and $S_{k+1}$ lie in the union $(\alpha_1, 0)\cup (\nu, \alpha_3)$.  Moreover, if consecutive zeros of  $S_i$ bracket the interval $[0, \nu]$, then there is a zero of the derivative of $S_i$ between $\nu$ and the larger of the two zeros.

If $S$ is a Stieltjes polynomial of degree $j$, then substitution into \eqref{lame1} and identification of powers of $x$ implies that
\begin{equation}
\mu = \mu_j = j\left( j-1 + 2(\rho_1+\rho_2+\rho_3)\right).
\label{mu}
\end{equation}
Thus $S_k$ and $S_{k+1}$ satisfy
\begin{eqnarray}
A(x)S_k'' + B(x)S_k' & = & \mu_k (x-\nu) S_k \label{k}\\
A(x)S_{k+1}'' + B(x)S_{k+1}' & = & \mu_{k+1} (x-\nu) S_{k+1} \label{k+1}
\end{eqnarray}
Define the integrating factor
\[ J(x) = \prod_{i=1}^3 \left| x - \alpha_i\right|^{2\rho_i}.
\]
Then $J' = J \frac{B}{A}$.  We derive two expressions, the first of which is obtained by multiplying \eqref{k} by $\mu_{k+1}S_{k+1}$ and \eqref{k+1} by $\mu_kS_k$ and taking the difference.
  The second is obtained by dividing the equations (\ref{k} - \ref{k+1}) by $A$, multiplying \eqref{k} by $S_{k+1}$, \eqref{k+1} by $S_k$, and taking the difference of the result.  The result is the following:
\begin{eqnarray}
\frac{d}{dx} \left[ J \left(\mu_{k+1}S_{k}' S_{k+1} - \mu_kS_{k} S_{k+1} ' \right)\right] & = & \left(\mu_{k+1} - \mu_k\right) J S_k' S_{k+1}',
\label{sturm2} \\
\mbox{and \ } \;\;\;\;\; \frac{d}{dx} \left[ J \left(S_{k+1}' S_k - S_{k+1} S_k ' \right)\right] & = & \left(\mu_{k+1} - \mu_k\right) Q S_k S_{k+1},
\label{sturm}
\end{eqnarray}
where $Q(x) = (x-\nu) J(x)/A(x)$.  Moreover, at any of the singular points $\alpha_1, 0$, or $\alpha_3$,
\begin{equation}
\mu_{k+1}S_{k}'S_{k+1} -  \mu_kS_{k}S_{k+1}' = 0.
\label{psizero}
\end{equation}

Notice that $Q(x) < 0$ for all $x\in (\alpha_1, 0)\cup(\nu, \alpha_3)$.  Now, consider two consecutive zeros of $S_k$, $x_1<x_2$ in either $(\alpha_1, 0)$ or $(\nu, \alpha_3)$.   Suppose that $S_k$ and
$S_{k+1}$ are positive in $(x_1, x_2)$.  Then
\[ \left.J\left(S_{k+1}'S_k-S_{k+1}S_k'\right)\right|_{x=x_1} \leq 0 \;\; \mbox{and} \;\;
\left.J\left(S_{k+1}'S_k-S_{k+1}S_k'\right)\right|_{x=x_2} \geq 0,
\]
but, according to \eqref{sturm}, $J\left(S_{k+1}'S_k-S_{k+1}S_k'\right)$ must be strictly decreasing on $(x_1, x_2)$, which is a contradiction.  Thus, $S_{k+1}$ must change sign on $(x_1, x_2)$.
Since $J(\alpha_i)=0$, a similar argument shows that there is a zero of $S_{k+1}$ (i) between any two zeros of $S_k$, (ii) between $\alpha_1$ and the smallest zero of $S_k$ in $(\alpha_1, 0)$, (iii) between the largest zero of $S_k$ in $(\alpha_1, 0)$ and 0, and (iv) between the largest zero of $S_k$ in $(0, \alpha_3)$ and $\alpha_3$.

There are two cases to consider, depending on whether there is a zero of $S_k$ less than zero.

\textit{Case 1:  There is a zero of $S_k$ in $(\alpha_1, 0)$.}  Firstly, if all the zeros of $S_k$ (and hence all the zeros of $S_{k+1}$) were to lie in $(\alpha_1, 0)$, then $JS_k'S_{k+1}'$ would not change sign in $(0, \alpha_3)$, which would contradict \eqref{sturm2}, since $J(0)=J(\alpha_3)=0$.

So, there must be at least one zero of $S_k$ in $(\nu, \alpha_3)$.  A zero of $S_k$ in $(\alpha_1, 0)$ implies that the zeros of $S_k$ and $S_{k+1}$ interlace, and hence that the zeros of the derivatives of $S_k$ and $S_{k+1}$ interlace \cite{rasc02}.
It is possible for a Van Vleck zero to equal zero, but since $A(0)=0$, this is possible if and only if the derivative of the corresponding Stieltjes polynomial has a zero at zero.  But, since the zeros of $S_k'$ and $S_{k+1}'$ interlace in this case, it is impossible for $\nu=0$ to be a Van Vleck zero corresponding to Stieltjes polynomials of successive orders.

Thus, it must be that $\nu>0$ and $S_k'(\nu), S_{k+1}'(\nu)\neq 0$.  And since the zeros of $S_k'$ and $S_{k+1}'$ interlace, the smallest zero of $S_k'$ in $(\nu, \alpha_3)$ is less than the smallest zero of $S_{k+1}'$ in this interval.  Let $\xi$ be the smallest zero of $S_k'$ in $(\nu, \alpha_3)$.  We may assume that $S_k, S_{k+1} > 0$ in $[0, \nu]$, so that $S_{k+1}'>0$ in $[0, \xi]$, and hence
\begin{equation}
 \left. \left(\mu_{k+1}S_{k}'S_{k+1} -  \mu_kS_{k}S_{k+1}'\right)\right|_{x=\xi}  =  -\mu_k S_k(\xi) S_{k+1}'(\xi) < 0.
\end{equation}
But this, along with \eqref{psizero}, contradicts \eqref{sturm2} since $JS_k'S_{k+1}' > 0$ in $(0, \xi)$.

\textit{Case 2:  All the zeros of $S_k$ lie in $(\nu, \alpha_3)$.}  In this case, there can be at most one zero of $S_{k+1}$ in $(\alpha_1, 0)$.  If there is one such zero of $S_{k+1}$ in $(\alpha_1, 0)$, then there is a zero of $S_{k+1}'$ between $\nu$ and the smallest zero of $S_{k+1}$ in $(\nu, \alpha_3)$.  Thus, whether there is one or no zeros of $S_{k+1}$ in $(\alpha_1, 0)$,
 $JS_k'S_{k+1}'$ does not change sign in $(\alpha_1, 0)$.  But this contradicts \eqref{sturm2} since $J(\alpha_1)=J(0)=0$.
\end{proof}

Now we are in a position to prove our main result.

\begin{theo}
Let $k$ be a positive integer.  Write the Van Vleck zeros of order $k$  in increasing order as
 $\nu_1^{(k)} < \nu_2^{(k)} < \cdots < \nu_{k+1}^{(k)}$.  Then
\begin{equation}
\alpha_1 < \nu_{1}^{(k+1)} < \nu_1^{(k)} < \nu_{2}^{(k+1)} < \nu_2^{(k)} < \cdots < \nu_{k+1}^{(k)} < \nu_{k+2}^{(k+1)} < \alpha_3.
\label{interlacing}
\end{equation}
In other words, the Van Vleck zeros of order $k$ and $k+1$ interlace:  between any two Van Vleck zeros of order $k$ there is a Van Vleck zero of order $k+1$, and vice versa.
\label{maintheorem}
\end{theo}
\begin{proof}
Under the transformation $x\mapsto (x-\alpha_2)/(\alpha_2 - \alpha_1), \; \nu\mapsto (\nu-\alpha_2)/(\alpha_2 - \alpha_1)$ the polynomial coefficients $A(x)$ and $B(x)$ in the Lam\'e equation \eqref{lame1} take the form
\begin{eqnarray*}
A(x)& =& x(x+1)(x-\alpha ), \\
 B(x) &=&  2(\rho_1+\rho_2+\rho_3) x^2 + 2\left(\rho_2+\rho_3 - \alpha(\rho_1+\rho_2)\right) x - 2\alpha \rho_2 ,
\end{eqnarray*}
where $\alpha = (\alpha_3 - \alpha_2)/(\alpha_2 - \alpha_1)$.   Now, suppose $S_k(x) = \sum_{j=0}^k a_j x^j$ is a degree $k$ polynomial solution of \eqref{lame1}.  Substitution into \eqref{lame1} and identification of the powers of $x$ yields the following relation for $j=0, \dots, k$:
\begin{eqnarray}
(\mu_{j-1} - \mu_k) a_{j-1} + j \left[(1-\alpha) (j-1) +  2\left(\rho_2+\rho_3 - \alpha(\rho_1+\rho_2)\right)\right] a_j&& \nonumber \\
 - (j+1)\alpha \left( j +2 \rho_2\right) a_{j+1}  =  -\mu_k \nu a_j ,&&
 \label{coefficients}
\end{eqnarray}
where $\mu_k$ is as in \eqref{mu} and $a_{-1}=a_{k+1}=0$.
Therefore, the coefficients of $S_k$ and the Van Vleck zeros of order $k$ are the eigenvectors and eigenvalues, respectively, of a matrix $B^{(k)}$, whose coefficients are functions of $\rho_1, \rho_2, \rho_3$, and $\alpha$.

Since the eigenvalues of a matrix are continuous functions of its entries \cite{hojo85}, and in light of Lemma \ref{unequal}, it suffices to show that for some particular values of $\rho_1, \rho_2, \rho_3$, and $\alpha>0$, the eigenvalues of $B^{(k)}$ and $B^{(k+1)}$ interlace.  So let $\rho_1=\rho_2=\rho_3 = 1/2$.  Then $B^{(k)}$ is tridiagonal, with nonzero entries given by
\begin{equation}
b_{j,j-1}^{(k)} = 1 - \frac{j(j-2)}{k(k+2)}, \;\;\; b_{j,j}^{(k)}=(\alpha-1)\frac{j(j-1)}{k(k+2)}, \;\;\;
b_{j,j+1}^{(k)} = \alpha \frac{j^2}{k(k+2)}
\end{equation}
$B^{(k)}$ is thus a function of $\alpha$, which we write as
\[B^{(k)}(\alpha) = B^{(k)}(0) + \alpha A^{(k)}.
\]
Note that $B^{(k)}(0)$ is bidiagonal, so its eigenvalues are given by the diagonal entries $-j(j-1)/k(k+2)$.  Since
\[\frac{j(j-1)}{(k+1)(k+3)} < \frac{j(j-1)}{k(k+2)} < \frac{j(j+1)}{(k+1)(k+3)}, \;\;\; j=2,3,\dots, k+1,
\]
the eigenvalues of $B^{(k)}(0)$ and $B^{(k+1)}(0)$ interlace except for a common eigenvalue at zero.  Call those eigenvalues of $B^{(k)}$ and $B^{(k+1)}$ that are zero at $\alpha=0$, $\lambda_0^{(k)}$ and $\lambda_0^{(k+1)}$, respectively.  These are functions of $\alpha$, and $\lambda_0^{(k)}(\alpha) = \nu_{k+1}^{(k)}$ when $\alpha>0$.
Therefore, for small $\alpha>0$, the eigenvalues of $B^{(k)}$ and $B^{(k+1)}$ interlace as long as
\begin{equation}
\left.\frac{d}{d\alpha} \lambda_0^{(k)}\right|_{\alpha=0} < \left.\frac{d}{d\alpha} \lambda_0^{(k+1)} \right|_{\alpha=0}.
\end{equation}
The rate of change of the eigenvalues at $\alpha=0$ are
\[\left.\frac{d}{d\alpha} \lambda_0^{(k)}\right|_{\alpha=0} = \frac{\mathbf{w}^TA^{(k)}\mathbf{v}}{\mathbf{w}^T\mathbf{v}},
\]
where $\mathbf{w}^T$ and $\mathbf{v}$ are left and right eigenvectors of $B^{(k)}(0)$ associated with the zero eigenvalue, and similarly for $\lambda_0^{(k+1)}$.  Since the entries of $B^{(k)}(0)$ in the first row are zero, the left eigenvector associated with the zero eigenvalue is  $\mathbf{w} = (1,0,\dots,0)^T$.  A simple calculation shows that the ratio of the second and first entries of $\mathbf{v}$ satisfy $v_2/v_1 = -b_{21}^{(k)}(0)/b_{22}^{(k)}(0)$.  The only nonzero entry of $A^{(k)}$ in the first row is $a_{12}^{(k)}$, so
\begin{eqnarray}
\left.\frac{d}{d\alpha} \lambda_0^{(k)}\right|_{\alpha=0}& =& a_{12}^{(k)} \frac{v_2}{v_1} \;\; = \;\;
-a_{12}^{(k)} \frac{b_{21}^{(k)}(0)}{b_{22}^{(k)}(0)} \;\; = \;\; \frac{k+2}{k+3} \nonumber \\
& < & \frac{k+3}{k+4} \;\; =  \left.\frac{d}{d\alpha} \lambda_0^{(k+1)}\right|_{\alpha=0}
\end{eqnarray}
Therefore, for small $\alpha>0$ and $\rho_i = 1/2$, the eigenvalues of $B^{(k)}$ and $B^{(k+1)}$ interlace.  Thus, the Van Vleck zeros of order $k$ and $k+1$ must interlace for any fixed positive $\rho_i$'s and any $\alpha_1<\alpha_2<\alpha_3$.   For, if there were a set of $\rho_i$'s and $\alpha_i$'s for which the interlacing \eqref{interlacing} did not hold, the continuity of the eigenvalues of a matrix with respect to its entries and the intermediate value theorem would imply the existence of a set of $\rho_i, \, \alpha_i$ for which $B^{(k)}$ and $B^{(k+1)}$ had a common eigenvalue, contradicting Lemma \ref{unequal}.
\end{proof}

\section{The non-orthogonality of the spectral polynomials}

The construction above results in a six parameter family of matrices $\left\{B^{(k)}\right\}_{k=1}^{\infty}$ for which the eigenvalues for successive $k$'s interlace and are in the interval $(\alpha_1,\alpha_3)$.  Note that $B^{(k)}$ is not a submatrix of $B^{(k+1)}$, so this result is not a simple consequence of the Cauchy interlacing theorem.  Moreover, consider the spectral polynomials of $B^{(k)}$ formed from the  Van Vleck zeros of order $k$:
\begin{equation} p_{k+1}(x) = \prod_{i=1}^{k+1} \left( x - \nu_i^{(k)}\right).
\label{spectralpoly}
\end{equation}
Then $p_k$ is a polynomial of degree $k$ with simple zeros that interlace with the zeros of $p_{k+1}$.  In this
section we show that the family $\left\{p_k\right\}$ is not orthogonal with respect to any measure. The theorem
depends on the following technical lemma.

\begin{lemma}  For any distinct numbers $\alpha_1, \alpha_2, \alpha_3$ and any $\rho_1, \rho_2, \rho_3 > 0$,
\begin{eqnarray}
\sum_{i=1}^{k+1} \nu_i^{(k)} & = & \frac{\alpha_1+\alpha_2+\alpha_3}{3}k + C_1(\alpha_j, \rho_j) + {\cal O}\left(\frac{1}{k}\right), \;\;\; \mbox{and \ } \label{trace1} \\
\sum_{i=1}^{k+1} \left(\nu_i^{(k)}\right)^2 & = &
 \left[\frac{\alpha_1^2+\alpha_2^2+\alpha_3^2}{5}  + \frac{2}{15}\left( \alpha_1\alpha_2 + \alpha_2\alpha_3 + \alpha_1\alpha_3\right) \right] k \nonumber \\
 && \; + \; C_2(\alpha_j, \rho_j) + {\cal O}\left(\frac{1}{k}\right) \label{trace2}
\end{eqnarray}
where $C_1$ and $C_2$ depend only on the $\alpha_j$'s and $\rho_j$'s.
\label{trace}
\end{lemma}
\begin{proof}
We make the same transformation as in the proof of Theorem \ref{interlacing}.  From \eqref{coefficients} we see that the transformed Van Vleck zeros are the eigenvalues of the tridiagonal matrix $B^{(k)}$ whose non-zero entries are
\[ b_{j,j-1}^{(k)} = 1 - \frac{\mu_{j-2}}{\mu_k}, \;\;\; b_{j,j}^{(k)}=\frac{(j-1)\left((\alpha-1)(j-2) +g_1\right)}{\mu_k},
\;\;\;
b_{j,j+1}^{(k)} = \alpha \frac{j(j+2\rho_2)}{\mu_k},
\]
where $g_1 = -2(\rho_2+\rho_3-\alpha(\rho_1+\rho_2))$.  The first equality \eqref{trace1} follows from a calculation of the trace:
\begin{eqnarray*}
\mbox{tr} \left(B^{(k)}\right)  &=&
\frac{(k+1)k(k-1) (\alpha -1)/3 +k(k+1)g_1/2}{k\left(k-1+2(\rho_1+\rho_2+\rho_3)\right)} \\
&=& \frac{\alpha-1}{3} (k+1) +  \frac{\alpha-1}{3} (1-2(\rho_1+\rho_2+\rho_3)) + \frac{g_1}{2} + {\cal O}\left(\frac{1}{k}\right)
\end{eqnarray*}
Making the transformation back to the original variables, $\nu \mapsto \nu (\alpha_2 - \alpha_2) + \alpha_2$,
\begin{eqnarray*}
\sum_{i=1}^{k+1} \nu_i^{(k)} & = & (\alpha_2 - \alpha_1)\, \mbox{tr} \left(B^{(k)}\right) + \alpha_2 (k+1) \\
& = & \left[(\alpha_2 - \alpha_1) \frac{\alpha-1}{3} + \alpha_2\right] (k+1) \\
&&+ (\alpha_2 - \alpha_1)\left[\frac{(\alpha-1)(1-2(\rho_1+\rho_2+\rho_3))}{3} + \frac{g_1}{2} \right] + {\cal O}\left(\frac{1}{k}\right)  ,
\end{eqnarray*}
and since $(\alpha_2 - \alpha_1)(\alpha - 1)/3 + \alpha_2 = \left(\alpha_1+\alpha_2+\alpha_3\right)/3$, this establishes \eqref{trace1}.

To prove \eqref{trace2} we must compute the trace of $\left(B^{(k)}\right)^2$.  The diagonal terms of this matrix are given by
\begin{eqnarray*}
\left(B^{(k)2}\right)_{jj} & = & \left\{  \left(\mu_k-\mu_{j-2}\right)\alpha (j-1)(j-1+2\rho_2)\right. \\
&& + \left. (j-1)^2\left((\alpha-1)(j-2) +g_1\right)^2
+\left(\mu_k-\mu_{j-1}\right)\alpha j(j+2\rho_2)  \right\}/ \mu_k^2 \\
&=&\frac{\left[(\alpha-1)^2-2\alpha\right] j^4 + {\cal O}(j^3)}{k^4 + {\cal O}(k^3)}
+ \frac{2\alpha j^2 + {\cal O}(j)}{k^2 + {\cal O}(k)}
\end{eqnarray*}
Since the only terms that are not constant (w.r.t. $k$) or ${\cal O}(1/k)$ come from the $j^4/k^4$ and $j^2/k^2$ terms,
\[ \mbox{tr} \left(B^{(k)2}\right)  =  \left[\frac{(\alpha-1)^2}{5}+\frac{4\alpha}{15}\right]\, k
 + \mbox{const.} + {\cal O}\left(\frac{1}{k}\right).
\]
As before, we transform back to the original variables,
\begin{eqnarray*}
 \sum_{i=1}^{k+1}\left( \nu_i^{(k)}\right)^2 & = & (\alpha_2 - \alpha_1)^2\, \mbox{tr}   \left(B^{(k)2}\right)  + 2\alpha_2 (\alpha_2 - \alpha_1)\, \mbox{tr}   \left(B^{(k)}\right) + \alpha_2^2 (k+1) \\
 &=& \left\{ (\alpha_2 - \alpha_1)^2  \left[\frac{(\alpha-1)^2}{5}+\frac{4\alpha}{15}\right]
 + 2\alpha_2 (\alpha_2 - \alpha_1)\frac{\alpha-1}{3} + \alpha_2^2 \right\}\, k \\
 && +\, C_2(\alpha_j, \rho_j) + {\cal O}\left(\frac{1}{k}\right).
 \end{eqnarray*}
 A simplification of the quantity in braces above establishes \eqref{trace2}
\end{proof}

\begin{theo}
The spectral polynomials \eqref{spectralpoly} are not orthogonal with respect to any measure.
\end{theo}
\begin{proof}
Suppose that $\{p_k\}$ is orthogonal with respect to some measure.  Then the polynomials must satisfy a
three-term recurrence relation of the following form \cite{sz75}.   There exists sequences $\{a_n\}$ and
$\{b_n\}$, with $a_n \in \R$ and $b_n > 0$ such that
\begin{equation}
p_n(x) = \left(x-a_n\right) p_{n-1}(x) - b_n\, p_{n-2}(x).
\label{3term}
\end{equation}
We will show that if the polynomials defined by \eqref{spectralpoly} satisfy the relation \eqref{3term}, then $ a_n$ converges to a real number and $b_n$ converges to a positive number.  This implies, by a well known theorem, that the density of zeros of $p_n$ in the limit as $n\rightarrow\infty$ is described by a measure which is different from the asymptotic density of zeros of spectral polynomials calculated by Borcea and Shapiro \cite{bosh08}.  This contradiction will imply the truth of the theorem.

First, identifying the coefficients of $x^n$ in both sides of \eqref{3term}, yields the equation
\[a_n = \sum_i \nu_i^{(n-1)} - \sum_i \nu_i^{(n-2)}.
\]
Therefore, Lemma \ref{trace} implies
\begin{equation}
\lim_{n\rightarrow\infty} a_n = \frac{\alpha_1+\alpha_2+\alpha_3}{3}.
\label{An}
\end{equation}
Next, we consider the coefficients of $x^{n-1}$. We get
\begin{eqnarray}
b_n & = & \left( \sum_i\nu_i^{(n-1)} - \sum_i\nu_i^{(n-2)} \right) \sum_i\nu_i^{(n-2)}
+ \sum_{i<j} \nu_i^{(n-2)}\nu_j^{(n-2)} -   \sum_{i<j} \nu_i^{(n-1)}\nu_j^{(n-1)} \nonumber \\
& = & -\frac{1}{2} \left(\sum_i\nu_i^{(n-1)} - \sum_i\nu_i^{(n-2)} \right)^2 +
\frac{1}{2}\left[\sum_i\left(\nu_i^{(n-1)}\right)^2 - \sum_i\left(\nu_i^{(n-2)}\right)^2 \right] \nonumber \\
&=& -\frac{1}{2}a_n^2 + \frac{1}{2}\left[\sum_i\left(\nu_i^{(n-1)}\right)^2 - \sum_i\left(\nu_i^{(n-2)}\right)^2 \right] \label{Bn}
\end{eqnarray}
Utilizing the second equality in Lemma \ref{trace} and
combining the result with \eqref{An} and \eqref{Bn}, we find
\[ \lim_{n\rightarrow\infty} b_n  = \frac{2}{45} \left( \alpha_1^2+\alpha_2^2+\alpha_3^2 -
\alpha_1\alpha_2 - \alpha_2\alpha_3 - \alpha_1\alpha_3\right),
\]
which is positive for all real $\alpha_1<\alpha_2<\alpha_3$.

Now, since $a_n \rightarrow a \in \R$ and $b_n \rightarrow b \in (0,\infty)$, according to Theorem 5.3 of \cite{ne79}, the polynomials $p_n$ have the asymptotic zero distribution $\omega_{[\alpha, \beta]}$ with density
\begin{equation}
\frac{d\omega_{[\alpha, \beta]}(x)}{dx} =
\left\{ \begin{array}{ll}
\frac{1}{\pi \sqrt{(\beta - x)(x-\alpha)}}, & \; \mbox{ \ if \ } \; x \in (\alpha, \beta) \\
0 & \;\; \mbox{ \ elsewhere}, \end{array}\right.
\label{asyorth}
\end{equation}
where $\alpha = a-2/b$ and $\beta = a+2/b$.  However, in \cite{bosh08} it is shown that the spectral polynomials $p_n$ defined by \eqref{spectralpoly} have the asymptotic zero distribution given by a probability measure supported on $(\alpha_1, \alpha_3)$, with density $\rho_A(x)$ given by
\begin{equation}
\rho_A(x) = \left\{\begin{array}{ll}
\frac{1}{2\pi} \int_{\alpha_2}^{\alpha_3} \frac{ds}{\sqrt{ (\alpha_3 - s) (s - \alpha_2) (s - \alpha_1) (s-x)}}
& \; \mbox{ \ if \ } \; \alpha_1 < x < \alpha_2, \\[.2in]
\frac{1}{2\pi} \int_{\alpha_1}^{\alpha_2} \frac{ds}{\sqrt{ (\alpha_3 - s) ( \alpha_2-s) (s - \alpha_1) (x-s)}}
& \; \mbox{ \ if \ } \; \alpha_2 < x < \alpha_3.
\end{array}\right.
\label{vvasy}
\end{equation}
Since the limiting distributions in \eqref{asyorth} and \eqref{vvasy} are unequal, it cannot be that the
spectral polynomials obey the recurrence relation \eqref{3term}, and hence they are not orthogonal with respect
to any measure. \end{proof}

\section{The Lam\'e equation}

The most common form of the Lam\'e equation in the literature is
\begin{equation} \label{Lame Eq}
 \frac{d^2 \phi}{dx^2} + (n(n+1) k^2 \text{sn}^2 \, (x,k) -h) \phi=0
\end{equation}
where $\text{sn} \, (x,k)$ is the Jacobi elliptic function with modulus $k$, $0<k<1$. For fixed $k$ and $n$, we
say that $h$ is an eigenvalue if \eqref{Lame Eq} admits a nontrivial solution. It is well known that if we
assume $n$ to be a positive integer (as we do from now on), then \eqref{Lame Eq} has  exactly $2n+1$ distinct
eigenvalues. Furthermore, the corresponding eigenfunctions are the Lam\'e functions of the first kind:
\begin{equation}
  \phi^{\gamma_1,\gamma_2,\gamma_3}(x) = \text{sn}^{\gamma_1} \, (x,k) \ \text{cn}^{\gamma_2} \, (x,k) \
  \text{dn}^{\gamma_3} \, (x,k) \ P_m(\text{sn}^2 (x,k))
\end{equation}
where $\gamma_i \in \{0,1\}$, and $P_m$ is a polynomial of degree $m$ with $n=2m+|\gamma|$.

The set $\Lambda_n$ of eigenvalues can be divided into eight disjoint subsets according to the different values
of $\gamma_i$, namely
\begin{equation*}
\Lambda_n = \begin{cases}
\Lambda_n^{0,0,0} \cup  \Lambda_n^{1,1,0} \cup  \Lambda_n^{1,0,1} \cup  \Lambda_n^{0,1,1} & \text{ if } n \text{ is even}\\
\Lambda_n^{1,0,0} \cup  \Lambda_n^{1,0,0} \cup \Lambda_n^{0,0,1} \cup  \Lambda_n^{1,1,1} & \text{ if } n \text{
is odd}
\end{cases}
\end{equation*}
where $\Lambda_n^{\gamma_1,\gamma_2,\gamma_3}$ is the set of all eigenvalues $h$  having an eigenfunction of the
form $\phi^{\gamma_1,\gamma_2,\gamma_3}$. The cardinality of each subset is
\begin{align*}
 |\Lambda_n^{0,0,0}|=n/2+1, \ |\Lambda_n^{1,0,0}|=|\Lambda_n^{0,1,0}|=|\Lambda_n^{0,0,1}|=(n+1)/2,\\
 |\Lambda_n^{1,1,0}|=|\Lambda_n^{1,0,1}|=|\Lambda_n^{0,1,1}|=n/2, \ |\Lambda_n^{1,1,1}|=(n-1)/2.
\end{align*}
Note that for $n$ even
\[ |\Lambda_n| = |\Lambda_n^{0,0,0}| + |\Lambda_n^{1,1,0}| + |\Lambda_n^{1,0,1}| + |\Lambda_n^{0,1,1}|=2n+1 \]
and similarly for $n$ odd.

Using the theory of Hill's equation, Volkmer \cite{vo04} recently obtained various interlacing properties
satisfied by the eigenvalues for a fixed value of the parameter $n$. We will use Theorem 2.1 to give new
interlacing results when $n$ takes consecutive integer values.

First, we need to rewrite Lam\'e equation into its algebraic form. Making the substitution $x \mapsto \text{sn
}^2 \, (x,k)$, we get
\begin{equation} \label{Lame Eq Alg}
  \frac{d^2 \phi}{dx^2} +\frac{1}{2} \left(\frac{1}{x}+\frac{1}{x-1}+\frac{1}{x-k^{-2}} \right)
   \frac{d \phi}{dx}= \frac{n(n+1)x- k^{-2} h }{4x(x-1)(x-k^{-2})} \phi.
\end{equation}
Now, if we substitute the corresponding Lam\'e function
\[ \phi^{\gamma_1,\gamma_2,\gamma_3}(x)= |x|^{\gamma_1/2} |x-1|^{\gamma_2/2} |x-k^{-2}|^{\gamma_3/2} P_m(x) \]
into \eqref{Lame Eq Alg}, then one can easily verify that the polynomial $P_m$ satisfies the Heine-Stieltjes
equation
\begin{multline*}
  \frac{d^2 P_m}{dx^2} + \left(\frac{\gamma_1+1/2}{x}+\frac{\gamma_2+1/2}{x-1}+\frac{\gamma_3+1/2}{x-k^{-2}}
  \right) \frac{d P_m}{dx} \\
  = \frac{m(m+|\gamma|+1/2)(x- \lambda)}{x(x-1)(x-k^{-2})} P_m,
\end{multline*}
where $\lambda(\gamma,k)$ is related to the eigenvalue $h$ through the relation
\begin{equation} \label{eigenvalue}
 \lambda = \frac{k^{-2} h -  (1+k^2)\gamma_1 - \gamma_2 - k^2 \gamma_3 - 2 k^2 \gamma_1 \gamma_3 - 2
 \gamma_1 \gamma_2}{m(m+|\gamma|+1/2)}.
\end{equation}
The next result is then an immediate consequence of Theorem 2.1.

\begin{prop} For each $\gamma=(\gamma_1,\gamma_2,\gamma_3) \in \{0,1\}^3$, let $h_{j,n}^{\gamma}$ denote the ordered
eigenvalues of \eqref{Lame Eq} in $\Lambda_n^{\gamma}$, and let $\lambda_{j,n}^{\gamma}$ be defined by
\eqref{eigenvalue} with $h=h^{\gamma}_{j,n}$. We have
\[ 0 < \lambda_{j,n+2}^{\gamma} < \lambda_{j,n}^{\gamma} < \lambda_{j+1,n+2}^{\gamma} < k^{-2} \]
for all $j=1,...,|\Lambda_n^{\gamma}|$.

\end{prop}

\section{Remarks and a conjecture}

It is a well-known fact that the zeros of orthogonal polynomials interlace \cite{sz75}.  The family $\{p_k\}$ of
spectral polynomials \eqref{spectralpoly} is an example of a family of polynomials with interlacing zeros, but
is not orthogonal, and is thus a counterexample  to the converse of the statement that the zeros of orthogonal
polynomials interlace.  We are unaware of any other example of such a family arising from a second order
differential equation.  This is all the more striking when we consider that the density function $\rho_A$
defined in \eqref{vvasy} satisfies, by a theorem in \cite{bosh08}, the following Heun differential equation:
\[ 8A(x) \rho_A''(x) + 8A'(x)\rho_A'(x) + A''(x)\rho_A(x) = 0,
\]
where $A(x) = (x-\alpha_1)(x-\alpha_2)(x-\alpha_3)$ is the function used to define the Lam\'e equation \eqref{lame1}.  It would be interesting to determine if there are any other spectral polynomials of second order differential operators of this kind, i.e. with interlacing zeros, non-orthogonality and with asymptotic density satisfying a Heun equation.

In \cite{mcboag08} we considered the Lam\'e equation in the case when $\alpha_1,\alpha_2,\alpha_3$ are the vertices of an equilateral triangle in the complex plane.  Since the Lam\'e equation is invariant under complex affine transformations, we may assume in this case that the $\alpha_i$'s are the third roots of unity, $\alpha_j=\exp\left( i(j-1) 2\pi/3\right)$, $j=1,2,3$.  In the special case when $\rho_1=\rho_2=\rho_3$, we found that the Van Vleck zeros of order $3k-1$ are of the form
\[
\lambda_{n}\, e^{\frac{2\pi i}{3}j}, \;\;\; n=1,\dots, k, \;\; j=0,1,2,
\]
where $\lambda_n \in (0,1)$ is real.  In other words, the Van Vleck zeros lie on the lines connecting the triangle incenter to its vertices.  When the order is $3k$ or $3k+1$, there is an additional Van Vleck zero at the center of the triangle.  Numerical evidence suggests an analog to Theorem \ref{maintheorem} in the complex case.  Let $\lambda_n^{(3k-1)}$, $n=1,\dots, k$ be the distance of the Van Vleck zeros of order $3k-1$ from the triangle incenter.  These are distinct.  If we label them in increasing order as
$\lambda_1^{(3k-1)} < \lambda_2^{(3k-1)}<\cdots < \lambda_k^{(3k-1)}$, we conjecture that these distances interlace with those of order $3k+2$:
\[0<
\lambda_1^{(3k+2)} < \lambda_1^{(3k-1)} < \lambda_2^{(3k+2)} < \lambda_2^{(3k-1)}<\cdots < \lambda_k^{(3k+2)} < \lambda_k^{(3k-1)} < \lambda_{k+1}^{(3k+2)} < 1.
\]
We note that, unlike in the real case this property only holds in the complex case when the $\rho_i$'s are all equal, since when the $\rho_i$'s are not all equal, the symmetry breaks and the Van Vleck zeros do not lie on the lines connecting the incenter with the vertices.

%\vskip12pt \noindent {\bf Acknowledgements:}

\bibliography{biblio}

\end{document}